\documentclass[leqno]{article}
\usepackage{amsmath,amssymb,a4wide,color}%,MnSymbol}
\usepackage{amsthm}
\usepackage[curve,matrix,arrow,cmtip]{xy}

\NoComputerModernTips

\newdir^{ (}{{}*!/-3pt/\dir^{(}}
\newdir_{ (}{{}*!/-3pt/\dir_{(}}

\def\OM{\mathchoice
  {\rlap{\kern3.2pt$\overline{\phantom{L}}$}M}
  {\rlap{\kern3.2pt$\overline{\phantom{L}}$}M}
  {\rlap{\kern2.4pt$\scriptstyle\overline{\phantom{L}}$}M}
  {\rlap{\kern1.8pt$\scriptscriptstyle\overline{\phantom{L}}$}M}}

\let\le\leqslant
\let\ge\geqslant

 \def\deg{\mathop{\rm deg}\nolimits}
\def\diag{\mathop{\rm diag}\nolimits}
\def\GL{\mathop{\rm GL}\nolimits}

\def\id{{\rm id}}

\let\phi\varphi
\let\theta\vartheta
\let\epsilon\varepsilon
\let\setminus\smallsetminus

\newtheorem{Thm}{Theorem}[section]
\newtheorem{Prop}[Thm]{Proposition}
\newtheorem{Lem}[Thm]{Lemma}
\newtheorem{Cor}[Thm]{Corollary}

\theoremstyle{definition}
\newtheorem{Def}[Thm]{Definition}
\newtheorem{Rem}[Thm]{Remark}
\newtheorem{Ex}[Thm]{Example}

\def\UseTheoremCounterForNextEquation{\setcounter{equation}{\value{Thm}}\addtocounter{Thm}{1}}

\def\qed{{\hskip0pt\unskip\unskip\nobreak\hfil\penalty50
          \hskip1em\hbox{}\nobreak\hfil
           {$\square$}
          \parfillskip=0pt\finalhyphendemerits=0
          \par}\medskip}

\newenvironment{Proof}
               {\noindent{\bf Proof.}\ }
               {\qed}

               {\noindent{\bf Proof of #1.}\ }
               {\qed}

\newcommand{\BC}{{\mathbb{C}}}
\newcommand{\BF}{{\mathbb{F}}}
\newcommand{\BP}{{\mathbb{P}}}
\newcommand{\BQ}{{\mathbb{Q}}}
\newcommand{\BR}{{\mathbb{R}}}
\newcommand{\BZ}{{\mathbb{Z}}}
\def\Cinf{{\BC}_\infty}
\def\Finf{F_\infty}

\newcommand{\sumprm}{\sideset{}{'}\sum}
\newcommand{\prodprm}{\sideset{}{'}\prod}

\newcommand{\qBinom}[2]{{\genfrac{[}{]}{0pt}{}{#1}{#2}}_{q^d}}
\begin{document}

\title{Hecke Operators for higher rank Drinfeld Modular Forms}

\author{Dirk Basson}

%\date{}
\maketitle

\begin{abstract}
We calculate the effect of simple Hecke operators on $u$-expansions of higher rank Drinfeld modular forms, the eigenvalue for the Drinfeld discriminant function $\Delta_t$ and show that a certain natural class of Hecke operators is completely multiplicative.
\end{abstract}

%{\advance\baselineskip by -6pt
%\tableofcontents
%}

%%%%%%%%%%%%%%%%%%%%%%%%%%%%%%%%%%%%%%%%%%%%%%%%%%%%%%%%%%%%%%%%%%%%%%%%%%%%%%%%%%%%%%%%%%%%%%%%%%%%%%%%%%%%%%%%%%%%%%%%%%%%%%%%%%%

%\newpage
\addtocounter{section}{-1}
\section{Introduction}
\label{Sec:Intro}

Hecke operators are an essential tool in the study of (classical) modular forms. The idea is roughly to think of a modular form $f$ to be defined on the set of all lattices in $\mathbb{C}$, in which case the Hecke operator $T_n$ ($n\in\mathbb{Z}$) returns a new function $T_nf$ which is the sum of the values of $f$ applied to all lattices of index $n$ inside the original lattice. This defines a linear operator on the vector space of modular forms that turns out to be a very valuable tool. The modern approach is to define these operators using double cosets, making it possible to generalize to many other contexts. It is in some sense enough to study the Hecke operators $T_p$, where $p$ is a prime number, because they are multiplicative ($T_{mn}=T_mT_n$ if $\gcd(m,n)=1$) and there is the recurrence relation (\cite[(5.10)]{DiamondShurman})
\[ T_{p^r} = T_pT_{p^{r-1}} - p^{k-1}\langle p\rangle T_{p^{r-2}}.\]

\medskip

Goss \cite{GossES} showed that one can make the same definition in a positive characteristic setting, starting with the polynomial ring over a finite field $A=\mathbb{F}_q[t]$ instead of the integers $\mathbb{Z}$. One can then define Drinfeld modular forms in this setting and Hecke operators in either of the two ways described above. Gekeler \cite{GekelerDMC} and \cite{GekelerCoeff} also carried out some computations of these Hecke operators. At this point only Drinfeld modular forms of rank 2 were considered.

%, embed $A$ inside its fraction field $F$ and then inside the completion $\Finf$ at the infinite place and finally inside the completion of an algebraic closure $\Cinf$ of $\Finf$. 

%\medskip

Recently, Pink \cite{Pink} constructed a Satake compactification of Drinfeld moduli spaces of dimension greater than 1, and this allowed him to define modular forms algebraically as global sections of certain ample line bundles. An analytic definition was given by Basson, Breuer, Pink \cite{BBP} where it is also proved, that these definitions agree (\cite[Thm 10.9]{BBP}). Hecke operators are also defined and the algebraic and analytic Hecke operators are compared (\cite[Thm 12.13]{BBP}). The purpose of this work is to provide one step further in the study of Hecke operators on higher rank Drinfeld modular forms.

%In \cite{BBP}, some progress is made toward understanding the coefficients of the $u$-expansions of these modular forms. One of the main difficulties in the rank $r>2$ theory is that the modular forms are functions in $r$ variables and the coefficients in this expansion are no longer constants, but functions of $r-1$ variables. Luckily they turn out to be (weak) Drinfeld modular forms of rank one less \cite[Thm 5.9]{BBP}. 

%In \cite[Ch 12]{BBP} the foundation is laid to start the study of Hecke operators in this setting. The purpose of this work is to provide one step forward in their study.

We would like to express arbitrary Hecke operators in terms of a few key operators, like in the classical case where it was reduced to $T_p$ where $p$ is a prime number. The higher rank case does provide an extra difficulty. In rank 2, a lattice of index $p^2$ inside another must either be $p$ times the original or form a quotient isomorphic to $\mathbb{Z}/p^2\mathbb{Z}$, but in rank 3 one has to keep track of other index-types as well. One approach to do this is followed by Shimura \cite{ShimAut} and we follow his strategy closely in section \ref{sct:HeckeAlgebra}. 

Let $p\in A$ be an irreducible polynomial. We shall use the language of double cosets later, but let us give the lattice theoretic interpretation here. When studying a Hecke operator in terms of lattices, for an operator of index $p^n$, we could either restrict to one index-type or add over all sublattices of index $p^n$. In the first, we would choose an index-type, say $(A/p^2)\times(A/p^2)\times(A/p)$, and denote the corresponding Hecke operator by $T(p^2,p^2,p,1,\ldots, 1)$, while if we sum over all the Hecke operators of index $p^5$, we use the notation $\mathcal{T}_{p^5}$. The composite of two operators like $T(p^2,p^2,p,1\ldots, 1)$ and $T(p^2,p,1,\ldots, 1)$ naturally breaks up into a sum of various different operators of the same determinant, in this example $p^8$. The operator $\mathcal{T}_p$ brings some symmetry which is very useful. It is also what corresponds most closely to the classical case and we shall focus most of our attention on it. 

\medskip

The main results of the paper are as follows: In section \ref{Sec:Computations} we calculate a set of representatives for the double coset and use it to calculate the effect of the Hecke operator on the $u$-expansion of a Drinfeld modular form (Theorem \ref{Thm:HeckeExp}). We note some cases where the calculations simplify and use it to prove that Hecke operators preserve the spaces of modular forms, cusp forms and double cusp forms, respectively. In section \ref{sec:eigen} we show that the rank $r$ discriminant function $\Delta$ is an eigenform with eigenvalue $p^{q^{r-2}(q-1)}$, generalizing Gekeler's result \cite{GekelerCoeff} and deduce that each coefficient form (associated to $\varphi_t$) is an eigenform too. Finally, in section \ref{sct:HeckeAlgebra}, we study the algebra of double coset operators, culminating in Theorem \ref{thm:HeckeMult} stating that the Hecke operators $\mathcal{T}_N$ act completely multiplicatively.

%\medskip

%For elliptic modular forms and for rank 2 Drinfeld modular forms there are simple formulae that can be used as the definition of a Hecke operator. In both cases this formula is obtained from a double coset operator by computing a set of representatives. In this work we start by defining Hecke operators as double coset operators and then calculate a concrete formula that can be used to determine some eigenforms. One simplification present in the characteristic $p$ setting is that the Hecke operators are completely multiplicative.

%It turns out that many of the most natural examples of higher rank Drinfeld modular forms are eigenforms. Most notably, the Eisenstein series of weight $k$ for $\GL_r(A)$ are eigenforms for $T_p$ ($p\in A$) with eigenvalue $p^k$, while the coefficient forms are all eigenforms for $T_p$ (acting on the spaces of appropriate weight) with eigenvalue $p^{q-1}$. This presents another contrast with the classical case where a modular form is determined by its eigenvalues alone. This property is already present in the rank 2 case, as explained by Gekeler \cite{GekelerCoeff}. 

%Introduce the contents of the article. 
%
%Some history, in particular the relationship to \cite{Pink} and \cite{BreuerPink}.
%
%
%A more specific introduction, including some definitions and forward references to results. Perhaps mention some obstacles that needed to be overcome?

%An outline of the contents of each section?
\subsection*{Acknowledgements}
I would like to thank Florian Breuer and Lenny Taelman for many useful discussions.

\section{Drinfeld modular forms of rank $r\ge 2$}
\label{sec:basics}
In this section we make the necessary recollections from \cite{BBP}. We restrict ourselves to the simplest case where $A$ is the polynomial ring and leave generalizations for the future.

Let $\BF_q$ be the finite field with $q$ elements and let $A=\BF_q[t]$. In analogy to the construction of the complex numbers from the integers ($\BZ\subset \BQ\subset \BR\subset \BC$), we let $F$ be the fraction field of $A$, we let $\Finf=\BF_q((t^{-1}))$, the completion of $F$ with respect to the absolute value at infinity, and we let $\Cinf$ be the completion of an algebraic closure of $\Finf$. We remark that the algebraic closure of $\Finf$ is of infinite degree over $F$, which is why a second completion is needed. It is also this fact that gives rise to the existence of modular forms of arbitrary rank --- there exist discrete $A$-submodules of $\Cinf$ of arbitrary rank.

We define the Drinfeld period domain $\Omega^r$ to be the complement of all the $\Finf$ hyperplanes in $\BP^{r-1}(\Cinf)$. It inherits the structure of a rigid analytic space from $\BP^{r-1}(\Cinf)$ as an admissible open subset. The sections of the global sheaf $\mathcal{O}_{\Omega^r}$ of this rigid analytic space are called the holomorphic functions on $\Omega^r$.

We shall write an element $\omega\in\Omega^r$ as a column vector, normalized so that $\omega=(\omega_1,\omega_2,\ldots,\omega_{r-1},1)^T$. 
%(Recall that $\period$ is the Carlitz period.)  
By definition of $\Omega^r$, the entries $\omega_1,\ldots, \omega_{r-1}$ are $\Finf$-linearly independent, and in particular non-zero.

There is an action of $\GL_r(F)$ on $\Omega^r$ given by 
\[\gamma(\omega) = j(\gamma,\omega)^{-1}\gamma\omega,\]
where $\gamma\omega$ denotes the matrix product, and $j(\gamma,\omega)$ is a factor used to normalize $\gamma\cdot\omega$ so that its last entry is 1. Explicitly, $j(\gamma,\omega)$ is 
the last entry of $\gamma\omega$.

Let $k$ 
%and $m$ be integers
be an integer and let $\gamma\in\GL_r(F)$. We define an operator on $\mathcal{O}_{\Omega^r}$ by
\[f|_{k}\gamma(\omega)=j(\gamma,\omega)^{-k}f(\gamma(\omega)).\]
It is easy to verify that this defines a right  action of $\GL_r(F)$ on $\mathcal{O}_{\Omega^r}$.

\begin{Def}
Let $\Gamma\subset\GL_r(A)$ be an arithmetic subgroup, and $k$ be an integer. 
A \emph{weak Drinfeld modular form of rank $r$ and weight $k$} is a holomorphic function $f\in\mathcal{O}_{\Omega^r}$ satisfying $f|_{k}\gamma(\omega)=f(\omega)$ for all $\omega\in\Omega^r$ and all $\gamma\in\Gamma$.
\end{Def}

For a weak modular form to be a modular form, it also needs to be bounded as $\omega$ tends to all the  ``cusps''. As in the classical case, this is defined by giving a series expansion for $f$ in terms of a parameter at infinity. Here follows a brief outline of the development in \cite[Ch. 4 \& 5]{BBP}:

Let $U(F)$ be the algebraic subgroup of $GL_r(F)$ consisting of matrices of the form

\begin{equation}\label{DefU}
\left(\begin{array}{c|c}
1 & *\ \ldots\ * \\ \hline
\begin{array}{c} 0 \\[-5pt] \vdots \\[-4pt] 0
\end{array}
& \id_{r-1}
\end{array}\right)
\end{equation}

and let $\Gamma_U=U(F)\cap \Gamma$. There is an isomorphism $\iota:F^{r-1}\to U(F)$, sending the $r-1$-tuple $v'\in F^{r-1}$ to the $r-1$ entries marked $*$ in \eqref{DefU}. Then $\Lambda'=\iota^{-1}(\Gamma_U)$ is an $\BF_p$-vector space, so we may form the exponential function 
\[e_{\Lambda'\omega'}(X)=X\prodprm_{\lambda\in\Lambda'\omega'}\left(1-\frac{X}{\lambda}\right)\]
and define the parameter
\[u_{\omega'}(\omega_1)=e_{\Lambda'\omega'}(\omega_1)^{-1}.\]
Note that since $\omega_1,\omega_2,\ldots, \omega_r$ are $\Finf$ linearly independent, there is never division by zero in the above expression.

\begin{Thm}[Basson, Breuer, Pink]
Let $f$ be a weak modular form of weight $k$ for the group $\Gamma<\GL_r(A)$. There exists a series expansion
\[f(\omega)=\sum_{n\in\BZ} f_n(\omega')u_{\omega'}(\omega_1)^n\]
that converges to $f$ on a neighbourhood of infinity (some suitable admissible open subset of $\Omega^r$; see \cite[Def. 4.12]{BBP}),
where each $f_n:\Omega^{r-1}\to\Cinf$ is a uniquely determined holomorphic function. 
\end{Thm}
\begin{Proof}\cite[Prop. 5.4]{BBP}\end{Proof}

\begin{Def}
A weak modular form $f$ is said to be \emph{holomorphic at infinity} if the functions $f_n$ in its $u$-expansion are identically zero for all $n<0$.

Let $\Gamma\subset\GL_r(F)$ be an arithmetic subgroup, and $k$ be an integer. 
A weak modular form  of rank $r$, weight $k$ is called a \emph{modular form $f$  of rank $r$ and weight $k$} if for every $\delta\in\GL_r(F)$, the function $f|_{k}\delta$ is holomorphic at infinity.
\end{Def}

Though the set $\GL_r(F)$ is infinite, it turns out that it suffices to check holomorphy at infinity for finitely many $\delta$. In particular, when $\Gamma=\GL_r(A)$ one needs to check it only when $\delta$ is the identity.

\medskip

In this work, we shall work mainly with the example of coefficient forms. The other main source of examples is the Eisenstein series, but we have nothing more to say than was already said in \cite[Chapter 14]{BBP}.

%We shall deal with three main examples of modular forms in this work: Eisenstein series, coefficient forms and exponential coefficient forms. These are all treated in \cite{BBP}, but it is our opinion that the reader will benefit from a slightly more concrete feel that comes from treating them in our special case.

%============FROM HERE==============
%
%In the following theorem the functions $g_i$ are called coefficient forms and will be defined in Section \ref{Sec:Computations}.
%\begin{Thm} \label{Thm:Generators}
%The graded ring of analytic modular forms for $\GL_r(A)$ is generated by $r$ algebraically independent functions $g_1, g_2,\ldots,g_r$ of weight $q-1$, $q^2-1$, \ldots, $q^r-1$, respectively.
%\end{Thm}
%\begin{Proof}\cite[Thm. 17.5(a)]{BBP}.\end{Proof}

%%%%%%%%%%%%%%%%%%%%%%%%%%%%%%%%%%%%%%%%%%%%%%%%%%%%%%%%%%%%%%%%%%%%%%%%%%%%%%%%%%%%%%%%%%%%%%%%%%%%%%%%%%%%%%%%%%%%%%%%%%%%%%%%%%%%%%%%%%%%%%%%%%%%%%%%%%%%%%%%%%%%%%%%%%%%%%%%%%%%%%%%%%%%%%%%%%%%%%%%%%%%%%

\section{Preliminary computations of $u$-expansions}\label{Sec:Computations}
%Let $N\in A$ be a non-unit and let $\Gamma(N)$ be the kernel of the natural projection $GL_r(A)\to \GL_r(A/NA)$. Let $[v]\in (N^{-1}/A)^r$. 
%In \cite{BBP}, Eisenstein series of weight $k$ for $\Gamma(N)$ were defined by
%\UseTheoremCounterForNextEquation
%\begin{equation}
 %E_{[v]}^k(\omega)=\sum_{\lambda\equiv v\pmod{A^r}} \lambda^{-k}.
%\end{equation}
%
%
%\begin{Prop}\label{prop:EisensteinExpansion}
%Let $[v]\in(N^{-1}/A)^r$ and let $[v']\in (N^{-1}/A)^{r-1}$ be its projection onto the last $r-1$ entries. Then the constant term in the $u$-expansion is
%\[E^k_{[v]}=\begin{cases}
%E^k_{[v']}+\textup{ higher terms.}&\text{if }v_1\in A\\
%0+\textup{ higher terms.}&\text{if }v_1\notin A.
%\end{cases}\]
%\end{Prop}
%\begin{Proof}
%This is a simplified form of \cite[Prop 13.10]{BBP}. Note that the $u$-parameters may be different, as the group used in \cite{BBP} to compute the expansion was the stabiliser of the lattice coset $\Lambda + v$, which is larger than the group $\Gamma(N)$ we use here. However, \cite[Prop 5.14]{BBP} assures us that the constant coefficients are not affected by this change.
%\end{Proof}

%From the Eisenstein series, \cite{BBP} constructed \emph{coefficient forms} as follows: 

Let $L\subset \Cinf$ be any $\mathbb{F}_q$-vector space. Let $\alpha_i$ be the coefficients of the exponential function associated to $L$ and $\beta_i$ is defined as the coefficient of $X^{q^i}$ in the power series expansion of the logarithm (compositional inverse of the exponential function) associated to $L$, i.e.
\[e_{L}(X)=X+\alpha_1X^q+\cdots;\qquad \log_{L}(X)=X+\beta_1X^q+\cdots.\]

\begin{Prop}\label{Prop:GossPolynomials}
Let $L\subset \Cinf$ be an $\mathbb{F}_q$-vector space.
\begin{enumerate}
	\item\label{GP:Inverse} As functions of $z\in\Cinf$ we have
	\[e_{L}(z)^{-1}=\sum_{\lambda\in L}(z+\lambda)^{-1}.\]
	\item\label{GP:Exist} There exists a series of polynomials $G_{k}$ depending\footnote{We refrain from using $L$ in the subscript here, but will use it later, when the context is less clear.} on $L$ such that
	\[\sum_{\lambda\in L}(z+\lambda)^{-k}=G_{k}(e_{L}(z)^{-1}).\]
	\item\label{GP:Recurrence} $G_k(X)=X(G_{k-1}(X)+\alpha_1G_{k-q}(X)+\alpha_2G_{k-q^2}(X)+\cdots)$.
	\item\label{GP:Monic} $G_k(X)$ is monic of degree $k$.
	\item\label{GP:PowerP} $G_{pk}(X)=G_k(X)^p$, where $p$ is the characteristic.
	\item\label{GP:SmallPowers} If $k\le q$, then $G_k(X)=X^k$.
	\item\label{GP:NoXTerm} $G_k(0)=0$, and if $k\ge 2$, then $X^2$ divides $G_k(X)$.
	\item\label{GP:Derivative} $X^2G'_k(X)=kG_{k+1}(X)$.
	\item\label{GP:Congruence} Each non-zero term in $G_k(X)$ has exponent congruent to $k\pmod{q-1}$.
	\item\label{GP:q^r-1} If $k=q^m-1$, then $G_k(X)=\sum_{0\le i<m} \beta_i X^{q^m-q^i}$.
\end{enumerate}
\end{Prop}
\begin{Proof}
Goss \cite[Prop.\ 6.6]{GossAlg}, and Gekeler \cite[Thm.\ 2.6]{GekelerZeroes}.
\end{Proof}

%\begin{Rem}
%The polynomials $G_k$ depend on the lattice $\Lambda'\omega'$ as can be seen from the recurrence relation in  (\ref{GP:Recurrence}) above, because the $\alpha_i$ are the coefficients of the exponential function associated to $\Lambda'\omega'$. If there is uncertainty about the lattice, we shall emphasise it by writing $G_{k,\Lambda'\omega'}$.
%\end{Rem}

%There is also a definition of an Eisenstein series for $\GL_r(A)$, where the sum is taken over all the elements in a lattice.

\begin{Def}
The Eisenstein series of weight $k$ for $\GL_r(A)$ is the function $E^k:\Omega^r\to\Cinf$
\[E^k(\omega)=\sumprm_{\lambda\in \omega\Lambda} \lambda^{-k}.\]
\end{Def}
This series converges whenever $k\ge 1$, and defines a non-zero function whenever $k$ is divisible by $q-1$. 

For any $\omega\in\Omega^r$, consider the Drinfeld module $\varphi$ associated to the lattice $\omega\cdot A^r$. We have
%For each $a\in A$ one has the product formula
%\begin{equation}\label{eq:DrinfeldModuleEquation}
%\varphi_a(X)=aX\prodprm_{v\in(N^{-1}A/A)^r} \left(1-E^1_v(\omega) X\right),
%\end{equation} 
%where the expression on the left is a polynomial of the form
\[\varphi_a(X)=aX+g_1(a,\omega)X^q+\cdots+g_{r\deg a}(a,\omega)X^{q^{r\deg a}}.\]
By varying $\omega$ one sees that each $g_i(a,\omega)$ becomes a holomorphic function of $\omega$ and turns out to be a Drinfeld modular form of weight $q^i-1$.
%that can be written as a symmetric polynomial in the Eisenstein series $E^1_v$ (for $v\in (N^{-1}A/A)^r$, $v\ne 0$). 
By \cite[Prop 15.12 (c)]{BBP} the first coefficient of the $u$-expansion of $g_i$ is the lower rank $g_i'$ when $i\le (r-1)\deg N$; and 0 if $i> (r-1)\deg N$.

We make one more definition to simplify the appearance of the $u$-expansion of $E^{k}$. For a non-zero $a\in A$, let $u_a$ be the $u$-expansion of $e_{\omega'\Lambda'}(a\omega_1)^{-1}$ when computed as follows:
\[e_{\omega'\Lambda'}(a\omega_1)=\varphi_a(e_{\omega'\Lambda'}(\omega))=\varphi_a(u^{-1}),\]
where $\varphi$ is the Drinfeld module associated to the lattice $\omega'\Lambda'$. This last expression is a polynomial in $u^{-1}$ which can also be written as some power of $u^{-1}$ times a polynomial in $u$. Also take out the leading coefficient $\Delta_a=g_{r\deg a}(a,\omega)$ so that the polynomial in $u$ has constant term 1. More precisely, let 
\[h_a(X)=\Delta_a^{-1}X^{q^{(r-1)\deg a}}\varphi_a(X^{-1}).\]
Then
\[\varphi_a(u^{-1})=\Delta_a u^{-q^{(r-1)\deg a}}h_a(u)\]
and, since $h_a(u)=1+\cdots$, we replace $u_a$ by its geometric series expansion to obtain
\begin{equation}
\label{eq:u_a}
u_a := e_{\omega'\Lambda'}(a\omega_1)^{-1}=\Delta_a^{-1}u^{q^{(r-1)\deg a}}h_a(u)^{-1}.
\end{equation}

\begin{Prop}\label{Prop:EisensteinSeriesExpansion}
The Eisenstein series $E^k(\omega)$ is a modular form for $\GL_r(A)$ and has the  $u$-expansion
\[E^k(\omega)=(E')^k(\omega')+\sum_{a\in A\setminus\{0\}} G_{k,A^{r-1}\omega'}(u_a) .\]
\end{Prop} 
\begin{Proof}
We have
\begin{align*}
E^k(\omega)&=\sum_{\substack{(a_1,\ldots,a_r)\in A^r\\ (a_1,\ldots,a_r)\ne (0,\ldots,0)}}(a_1\omega_1+a_2\omega_2+\cdots+a_r\omega_r)^{-k}\\
&=\sum_{\substack{(a_2,\ldots,a_r)\in A^{r-1}\\ (a_2,\ldots,a_r)\ne (0,\ldots,0)}}(a_2\omega_2+\cdots+a_r\omega_r)^{-k}+\sum_{a_1\in A\setminus\{0\}}\sum_{(a_2,\ldots,a_r)\in A^{r-1}}(a_1\omega_1+a_2\omega_2+\cdots+a_r\omega_r)^{-k}\\
&=(E')^k(\omega')+\sum_{a_1\in A\setminus\{0\}} G_{k,A^{r-1}\omega'}(e_{\omega'\Lambda'}(a_1\omega_1)^{-1})\\
&=(E')^k(\omega')+\sumprm_{a_1\in A\setminus\{0\}} G_{k,A^{r-1}\omega'}(u_{a_1}).
\end{align*}
\end{Proof}
\begin{Rem}
This expression looks very simple, but that  is because all the computations are hidden in the definition of $u_a$. 
\end{Rem}

Before the next Proposition, it is also useful to recall the relations (for all $a\in A$ and $k\ge 1$)
\begin{equation}
\label{eq:RelatingEisensteinCoeff}
\left(a-a^{q^k}\right)E^{q^k-1}(\omega)=\sum_{i=1}^{k-1} E^{q^i-1}(\omega)g_{k-i}(a,\omega)^{q^i} 
\end{equation}
between the Eisenstein series and the coefficient forms from  \cite[(2.10)]{GekelerCoeff}.

\begin{Prop}\label{Prop:CoefficientsCongruent}
Let $f$ be one of the following modular forms for $\GL_r(A)$: An Eisenstein series $E^{q^k-1}(\omega)$ or a coefficient form $g_i(a,\omega)$ for any $a\in A\setminus \BF_q$. Let $f$ have the $u$-expansion $f(\omega)=\sum_{n\ge 0}f_n(\omega')$. If $f_m(\omega')$ is not identically 0, then $q-1\mid m$ and $m\equiv -1,0\pmod{q}$.
\end{Prop}
\begin{Proof}
The first statement follows from the modular functional equation for the scalar matrix $\gamma=c\cdot I$, with $c\in\mathbb{F}_q^\times$, applied to $f$ and to each $f_n$, remembering that each $f_n$ is also a weak modular form.
%The first statement is Lemma \ref{lem:CoefficientDivisibleByq-1}. 
Define the statement (depending on the function $f$)
\begin{itemize}
	\item[($\ast$)] for all $m$, if $f_m$ is not identically 0, then $m\equiv-1,0\pmod{q}$.
\end{itemize}
Our strategy is to establish $(*)$ for various power series related to the Eisenstein series, similar to Gekeler's strategy in \cite{GekelerCoeff}.

The polynomial $h_a(u)$ viewed as a power series satisfies $(*)$, since its non-zero terms have exponents $q^{\deg a}-q^m$ for $m=0,\ldots,\deg a-1$. Then also $u_a^{q^k-1}$ satisfies $(*)$, since $u_a^{q^k-1}=u_a^{q^k}u^{-q^{(r-1)\deg a}}h_a(u)$.  By Proposition \ref{Prop:GossPolynomials}(\ref{GP:Congruence}), the Goss polynomial $P_{q^k-1}(X)$ has non-zero coefficients only for those $X$ whose exponent is divisible by $q$, with the exception of the leading term $X^{q^k-1}$. Thus $P_{q^k-1}(u_a)$ satisfies $(*)$. The expansion from Proposition \ref{Prop:EisensteinSeriesExpansion} now shows that each Eisenstein series $E^{q^k-1}(\omega)$ satisfies $(*)$, and the result for the coefficient forms follows from the relations  \eqref{eq:RelatingEisensteinCoeff}.
\end{Proof}

\section{Hecke operators} 
\label{sec:Hecke}
We define Hecke operators as double coset operators as in the classical case. The following definition is taken from \cite[Def. 12.11]{BBP}:

\begin{Def}
Let $\Gamma$ and $\Gamma'$ be congruence subgroups of $\GL_r(F)$ and let $\delta\in \GL_r(F)$. We define the Hecke operator associated to $\delta$ as
\[ T_{\delta} : \mathcal{M}(\Gamma') \to \mathcal{M}(\Gamma), \quad f \mapsto \sum_{\gamma} f|_k \gamma,\]
where $\gamma$ runs through a set of coset representatives for $\Gamma'\backslash \Gamma' \delta\Gamma$.
\end{Def}

%Hecke operators are defined as double coset operators on spaces of modular forms. Let $\alpha\in \GL_r(F)$ and suppose that $\Gamma_1$ and $\Gamma_2$ are congruence subgroups of $\GL_r(A)$. Then the set
%\[\Gamma_1\alpha\Gamma_2=\{\gamma_1\alpha\gamma_2:\gamma_1\in\Gamma_1, \gamma_2\in\Gamma_2\}\]
%is a double coset. It can be decomposed into a finite union of cosets $\Gamma_1\beta_i$, so that
%\[\Gamma_1\backslash \Gamma_1\alpha\Gamma_2=\bigcup_i \Gamma_1\beta_i.\]
%Set $\Gamma_3=\alpha^{-1}\Gamma_1\alpha\cap \Gamma_2$, which also turns out to be a congruence subgroup. 
%
%\begin{Def}
%We define the  $\Gamma_1\alpha\Gamma_2:\M_k(\Gamma_1)\to\M_k(\Gamma_2)$ operator of weight $k$ by
%\[f[\Gamma_1\alpha\Gamma_2]_k=\sum_i f|_{k}\beta_i,\]
%where $\Gamma_1\backslash \Gamma_1\alpha\Gamma_2=\bigcup_i \Gamma_1\beta_i$.
%\end{Def}
%It is routine to check that the sum on the right is $\Gamma_2$ invariant, and holomorphy at infinity follows from the holomorphy at infinity of $(f|_k\beta_i)|_k\delta$ for every $\delta\in\GL_r(F)$. \\

%Because of the complexity of the extra dimensions in the theory of higher rank Drinfeld modular forms, 
We shall restrict ourselves to the simplest case where $\Gamma_1=\Gamma_2=\GL_r(A)$ and (with the exception of section \ref{sct:HeckeAlgebra})
\UseTheoremCounterForNextEquation
\begin{equation}
\label{eq:DeltaForm}
 \delta = \begin{pmatrix} p&0\\0&\id_{r-1}\end{pmatrix}
\end{equation}
is a diagonal matrix with 1's everywhere except in the top left entry which is some irreducible $p\in A$. We now describe a set of coset representatives used to compute this Hecke operator.
%As in the classical case, this corresponds to a more elementary definition of a Hecke operator which we describe by handling the double coset representatives by turn.

For any $m\in\{1,2,\ldots, r\}$ and column vector $\underline{b} = (b_1,b_2,\ldots, b_r)^T$, define the $r\times r$ matrix $\beta_{m,\underline{b}}$ to be the $r\times r$ identity matrix, except for column $m$, which is set equal to the column vector $\underline{b}$. For any $m\in\{1,2,\ldots, r\}$, let $B_m$ be the set
\[ B_m := \{\beta_{m,\underline{(b_1,\ldots, b_r)}}\mid b_{m+1}=\cdots = b_{r}=0,\ b_m = p,\ \forall i<m, \deg b_i<\deg p\}\]
be the set of $r\times r$ matrices $\beta_{m,\underline{b}}$ where the only non-trivial column $m$ has the entry $p$ on the diagonal, has zeroes below the diagonal, and the elements above the diagonal have degree less than the degree of $p$. Thus, a typical element in $B_m$ is an upper triangular matrix of the form
\[ \begin{pmatrix}
1&0&\cdots&b_1&&0 \\
&1&\cdots&b_2&&0 \\
&&\ddots&\vdots&&\vdots\\
&&& p&&0 \\
&&&&\ddots&0\\
&&&&&1
\end{pmatrix}.\]

\begin{Prop}
The union $\displaystyle \cup_{m=1}^r B_m$ is a set of coset representatives for $\GL_r(A)\backslash \GL_r(A)\delta\GL_r(A)$. 
\end{Prop}
\begin{proof}
First, note that each such matrix $\beta_{m,\underline{b}}$ is in the double coset $\GL_r(A)\delta\GL_r(A)$ by writing it as the product
\[ \beta_{m,\underline{b}} = P\delta (P \alpha_{m,\underline{b}}),\]
where $P$ is the $r\times r$ permutation matrix that switches rows 1 and $m$ and $\alpha_{i,\underline{b}}$ is the same as $\beta_{m,\underline{b}}$, except the $m$-th diagonal entry is 1 instead of $p$, and so $\alpha_{m,\underline{b}}$ (and thus $P\alpha_{m,\underline{b}}$) is in $\GL_r(A)$.

Secondly, we prove that no two of the $b_{m,\underline{b}}$ are in the same coset by showing that the row span (as $A$-modules) of these matrices are all different. Note that $m$ is uniquely determined as the first index from the end for which the $m$-th elementary unit row vector $e_m=(0,\ldots, 0, 1, 0, \ldots, 0)$ is not in the row span. Then, if $\beta_{m,\underline{b}}$ and $\beta_{m,\underline{c}}$ had the same row span, it would contain both the row vector $(0,\ldots, 0, 1, 0,\ldots, b_i, \ldots, 0)$ (with $1$ in column $i<m$ and $b_i$ in column $m$) and the row vector $(0,\ldots, 0, 1, 0,\ldots, c_i, \ldots, 0)$ and thus contains their difference $(0,\ldots, 0, b_i-c_i,\ldots 0)$. We have $\deg b_i<\deg p$ and $\deg c_i < p$, so if $b_i\ne c_i$, then $\deg (b_i-c_i)<\deg p$. But together with the $m$-th vector $(0,\ldots, 0,p,0,\ldots, 0)$ this would generate $(0,\ldots, 1, \ldots, 0)$, since $p$ is irreducible, contradicting the choice of $m$. Thus $b_i=c_i$ for all $i<m$, proving that $\beta_{m,\underline{b}} = \beta_{m,\underline{c}}$.

\medskip

Lastly, we prove that every element in the double coset is in one of these cosets. Note that every element in the double coset has entries in $A$ and has determinant a unit times $p$. We show that every matrix $\alpha$ with these properties is in some coset $\GL_r(A)\beta_{m,\underline{b}}$. Again we use the $A$-module spanned by the row vectors of $\alpha$. Since $\alpha$ is invertible, this module, call it $M$, is a submodule of $A^r$ of full rank, so each intersection $M_m:= M\cap(\{0\}^{m-1}\times A\times \{0\}^{r-m})$ (containing the elements whose only non-zero entry is in entry $m$) is non-zero. 

We claim that for each $m$, the module $M_m$ is either generated by $(0,\ldots, 0,1,\ldots, 0)$ or by $(0,\ldots, 0, p, \ldots, 0, \ldots, 0)$. If it were generated by some other $(0,\ldots, 0, p',0\ldots, 0)$, then $p'$ would have to divide the determinant of $\alpha$, which equals $p$, and this is impossible since $p$ is irreducible. Starting from the back, if $M_r$ is generated by $(0,0,\ldots, 1)$, that means that $(0,0,\ldots, 1)$ is in $M_r$ and thus there is some matrix in $\GL_r(A)$ that multiplied with $\alpha$ results in a last row of  $(0,0,\ldots, 1)$. Similarly, if $(0,\ldots,1,0)$ generates $M_{r-1}$, then $\alpha$ is row equivalent to a matrix with a second to last row of $(0,\ldots,1, 0)$. This continues until there is an $m$ where $M_m$ is generated by $(0,\ldots, p, \ldots, 0)$ instead of $(0,\ldots, 1,\ldots,0)$. 

Collect everything up to this point in a matrix $X\in\GL_r(A)$ such that $X\alpha$ has its last $r-m$ rows elementary unit vectors, row $m$ as $(0,\ldots, p, \ldots,0)$ and consider the top-left $(m-1)\times (m-1)$ submatrix $Y'$ of $X\alpha$. Its determinant is a unit in $A$, since the determinant of $\alpha$ equals $p$ times the determinant of the smaller matrix. That means it is an element of $\GL_{r-m}(A)$, and thus the matrix
\[ Y = \begin{pmatrix} Y&0\\ 0& \id_{r-m+1}\end{pmatrix}\]
is an element of $\GL_r(A)$ and $Y^{-1}X\alpha$ is the identity matrix, except for column $m$.

At this point the entries in column $m$ below the diagonal are already equal to 0, because of the way the previous rows were constructed, and the entries below the diagonal can be reduced modulo $p$ until the $b_i$ that remains has degree less than the degree of $p$.
\end{proof}

\begin{Rem}
It is also possible to find similar representatives that are lower triangular, but we choose upper triangular representatives, because then $\beta_{m,\underline{b}}\cdot\omega$ has the same last entry as $\omega$ (in most cases).
\end{Rem}

\begin{Ex}
In the rank 2 case, the representatives are
\[ \begin{pmatrix} p&0\\0&1\end{pmatrix},\quad \begin{pmatrix} 1&b\\ 0&p\end{pmatrix}\]
yielding the formula (using $\omega_2=1$)
\[ T_\delta f(\omega_1) = f(p\omega_1)+p^{-k}\sum_{\deg b<\deg p}f\left(\frac{\omega_1+b}{p}\right),\]
which is a scalar multiple ($p^{-k}$ times) the one in Gekeler \cite[(7.1)]{GekelerCoeff}.
\end{Ex}

%\begin{Cor}\label{Prop:HeckeActionComputation}
%Let $p\in A$ be an irreducible element. The Hecke operator $T_p:\CM_k(\GL_r(A))\to\CM_k(\GL_r(A))$ may be computed as follows:
%\[T_pf(\omega)=\sum_{m<r,\ \beta_{m,\underline{b}}\in B_m} f(\beta_{m,\underline{b}}\omega) + \sum_{\beta_{r,\underline{b}}\in B_r} p^kf((p\omega_1,\ldots,p\omega_{r-1},1)).\]
%\end{Cor}
%\begin{Proof}
%The first summand is clear --- it corresponds to the coset representatives $\beta_{m,\underline{b}}$ where $m<r$. In that case the last entry of $\beta_{m,\underline{b}}$ is equal to that of $\omega$, and no further scaling is needed.
%
%The terms in the second summand correspond to the case where the bottom right entry of $\beta_{r,\underline{b}}$ equals $p$ and thus the last entry of the product \todo CHECK THIS AGAINST THE DEFINITION!!
%\end{Proof}

%Next, we would like to understand what the effect of this Hecke operator is on the $u$-expansion of a modular form. The coefficients in such a $u$-expansion are modular forms of rank $r-1$, so it is perhaps natural to expect an expression involving lower rank Hecke operators. We define $\delta'$ as the $(r-1)\times (r-1)$ diagonal matrix with $p$ in its top left entry and 1 in every other diagonal entry; the rank $r-1$ version of $\delta$. We indicate the weights of the Hecke operators in the formula below, since the coefficients $f_n$ are lower rank (weak) modular forms of weight $k-n$, instead of $n$.

\medskip

Next we aim to describe how the Hecke operator $T_\delta$ acts in relation with the $u$-expansion of the modular form. 

The parameter $u$ in the expansion of a modular form depends on the lattice $\Lambda'= A^{r-1}\omega'$. When we act on the $u$-expansion with a matrix $\beta_{m,b}$, this lattice changes. Thus, we define
\[ \Lambda'_b = A^{r-1}(\beta_{m,b}\omega').\]
In the course of proving the next Proposition, we shall use Goss polynomials and therefore we also need the following finite $\BF_q$-vector spaces. For $m=2,\ldots, r-1$ we define
\[ L_b = e_{\Lambda'_b}(\Lambda')\]
and for $m=r$, we define 
\[ L_b = e_{p\Lambda'_b}(\Lambda').\]
The reason for the difference is that if $m=r$, then $j(\gamma, \omega)=p$ so the action of $\beta_{m,b}$ includes scaling down by $p$. We point it out again during the proof.

\begin{Thm}
\label{Thm:HeckeExp}
Let $f\in\mathcal{M}_k(\GL_r(A))$ and suppose that $f(\omega)$ has the $u$-expansion
\[ f(\omega) = \sum_{n} f_n(\omega')u(\omega_1)^n.\]
Let $T_{\delta,k}$ be the weight $k$ Hecke operator used above. Then $T_{\delta,k} f$ has the $u$-expansion
\[ T_{\delta,k}f(\omega) = \sum_n f_n(\omega')u(p\omega_1)^n + \sum_n \sum_{\beta_{m,b}\in \cup_{m=2}^{r}B_m}f_n|_{k-n}\beta_{m,\underline{b}}(\omega')\cdot G_{n,L_b}(u(\omega_1)),\]
where $G_{n,L_b}(X)$ is the Goss polynomial associated to the lattice $L_b$.
\end{Thm}
\begin{Rem}
Note that $L_b$ is a finite-dimensional $\mathbb{F}_q$-vector space. In \cite[(3.9)]{GekelerCoeff}, Gekeler explains why the Goss polynomials $G_{n,L}$ associated to an $\mathbb{F}_q$-vector space $L$ of dimension $d$ are divisible by $X^m$ for $m=\lfloor n\cdot q^{-d}\rfloor+1$. This implies that only finitely many terms in the second summand contribute to any given coefficient $u^m$.
\end{Rem}
\begin{Proof}
The first summand comes from the term $f|_k \beta_{1,\underline{0}}(\omega)$, when $\beta_{1,\underline{0}}\omega$ is unchanged, except for the first entry, which is $p\omega_1$.

Next, let us fix $m$ and study the $u$-expansion of the sum $\sum_{\beta_{m, \underline{b}}\in B_m}f|_k \beta_{m,\underline{b}}$. First assume that $m<r$. 
Then $j(\gamma,\omega) = 1$, so $f|_k\beta_{m,\underline{b}}(\omega) = f(\beta_{m,\underline{b}}\omega)$ and that
\[ \beta_{m,(b_1,\ldots, b_{m-1},1,0,\ldots,0)}\omega = (\omega_1+b_1\omega_m,\  \omega_2+b_2\omega_2,\ \ldots,\ \omega_{m-1}+b_{m-1}\omega_m,\ p\omega_m,\ \omega_{m+1},\ \ldots,\ \omega_r)^T.\]
Since $u_{\omega'}(\omega_1)$ is the inverse of the exponential function associated to the lattice $\Lambda'$ generated by the last $r-1$ entries of $\omega$, the effect on the $u$-expansion will involve the exponential function associated to the lattice $\Lambda'_{\underline{b}}$ generated by the last $r-1$ entries of $\beta_{m,\underline{b}}$. Explicitly
\[ \Lambda'_{\underline{b}} = A(\omega_2+b_2\omega_m)+\cdots+A(\omega_{m-1}+b_{m-1}\omega_m)+A(p\omega_m)+A(\omega_{m+1})+\cdots + A(\omega_r),\]
from which it is clear that $\Lambda'_{\underline{b}}$ is a sublattice of $\Lambda'$ with index $p$, i.e. $\Lambda'/\Lambda'_{\underline{b}}\cong A/pA$.

Let us now break up $\underline{b}$ as $(b_1,b')$ and consider the terms where $b'$ is constant. Summing over all $b_1$ of degree less than $\deg p$, we get
\[
	\sum_{b_1}f_n((\beta_{m,\underline{b}}\omega)') u_{\Lambda'_b}(\omega_1+b_1\omega_m)^n\]
where the factor $f_n((\beta_{m,\underline{b}}\omega)')$ doesn't depend on $b_1$, so we calculate the other factor:
\UseTheoremCounterForNextEquation
\begin{align}
\label{eq:sum_of_u_power_n}
	\sum_{b_1} u_{\Lambda'_b}(\omega_1+b_1\omega_m)^n &= \sum_{b_1} e_{\Lambda'_b}(\omega_1+b_1\omega_m)^{-n}\\ \UseTheoremCounterForNextEquation
		\label{eq:pre_GossPoly}
		&= \sum_{b_1} (e_{\Lambda'_b}(\omega_1)+e_{\Lambda'_b}(b_1\omega_m))^{-n}. 
\end{align}
%As $b_1$ varies over all elements of degree less than $\deg p$, the quantities $e_{\Lambda'_b}(b_1)$ vary through the finite $\BF_q$-vector space $L_{b} = e_{\Lambda'_b}(b_1A)$.
To apply Proposition \ref{Prop:GossPolynomials} 
we first evaluate this expression in the case where $n=1$. Then
\[ \sum_{b_1} e_{\Lambda'_b}(\omega_1+b_1\omega_m)^{-1} = \sum_{b_1}\sum_{\lambda\in\Lambda'_b} (\omega_1+\lambda+b_1\omega_m)^{-1}. \]
Now, as $\lambda$ varies through $\Lambda'_b$ and $b_1$ varies through elements of $A$ of degree less than $\deg p$, the sum $b_1\omega_m+\lambda$ varies through the full lattice $\Lambda'$ and thus
\[ \sum_{b_1}\sum_{\lambda\in\Lambda'_b} (\omega_1+\lambda+b_1\omega_m)^{-1} = e_{\Lambda'}(\omega_1)^{-1}=u.\]

Now, applying Proposition \ref{Prop:GossPolynomials} to \eqref{eq:pre_GossPoly} yields $G_{n,L_b}(u)$, because $L_b = e_{\Lambda'_b}(b_1A)$.

%\[ \sum_{b_1}\sum_{\lambda\in\Lambda'_b} (\omega_1+b_1\omega_m+\lambda)^{-n} = G_{n,\Lambda'}(u).\]

\medskip

Now we treat the case $m=r$, which is similar, but has some differences because $j(\gamma,\omega) = p$, so 
\[ \beta_{m,\underline{b}}(\omega) = \left(\frac{\omega_1+b_1}{p}\ldots, \frac{\omega_{r-1}+b_{r-1}}{p},1\right)^T\]
after normalization. First note that in the definition of $T_{\delta,k}(f)$, the terms $f|_k \beta_{r,\underline{b}}(\omega) = j(\beta_{m,\underline{b}})^{-k}f(\beta_{r,\underline{b}}(\omega)')$, which means that we have to add a factor of $p^{-k}$ which wasn't there in the previous case.

The lattice 
\[\Lambda'_b = A\left(\frac{\omega_1+b_1}{p}\right)+\cdots +A\left(\frac{\omega_{r-1}+b_{r-1}}{p}\right)+A\]
this times contains $\Lambda'$, but it is a sublattice of index $p$ of $\frac{1}{p}\Lambda'$. Similar to the computation to the previous case, we sum over the terms where $b'$ is fixed. 
\begin{align*}
	 \sum_{b_1}p^{-k}f_n(\beta_{m,\underline{b}}(\omega)')u_{\beta_{m,\underline{b}}(\omega)'}\left(\frac{\omega_1+b_1}{p}\right)^n &= p^{-k}f_n(\beta_{m,\underline{b}}(\omega)')\sum_{b_1} e_{\Lambda'_b}\left(\frac{\omega_1+b_1}{p}\right)^{-n}\\
	&= p^{-k}f_n(\beta_{m,\underline{b}}(\omega)')\cdot\left(\frac{1}{p}\right)^{-n}\sum_{b_1} e_{p\Lambda'_b}(\omega_1+b_1)^{-n}\\
	&= p^{-k}f_n(\beta_{m,\underline{b}}(\omega)')p^nG_{n,L_b}(u),
\end{align*}
as in the previous case, because $p\Lambda'_b$ does have index $p$ in $\Lambda'$, so the same argument works. Also note that in this case $f_n(\beta_{m,\underline{b}}(\omega)') = j(\beta_{m,\underline{b}},\omega')^{k-n} f_n|_{k-n}\beta_{m,\underline{b}}\ (\omega')$, so
\[ p^{n-k}f_n(\beta_{m,\underline{b}}(\omega)')G_n(u) = f_n|_{k-n}\beta_{m,\underline{b}}\ (\omega')G_n(u).\]
%
%\medskip
%
%\begin{align*}
  %\sum_{2\le m\le r}\sum_{b} &f_n|_{k-n}\beta_{m,\underline{b}}\ (\omega')u_{\beta_{m,\underline{b}}(\omega)'}(\omega_1+b_1\omega_m)^n \\
	%&= \sum_{2\le m\le r}\sum_{b'}f_n|_{k-n}\beta_{m,\underline{b}}\ (\omega')G_{n,\Lambda'}(u) \\
	%&= T_{\delta',k-n}(f_n)(\omega')G_{n}(u).
%\end{align*}
\end{Proof}

\begin{Ex}
Again, this result specializes to Gekeler's formula \cite[(7.3)]{GekelerCoeff}:
\[ T_{\wp}\left(\sum a_it^i\right) = p^k\sum a_1t_p^i + \sum a_iG_{i,\wp}(pt).\] 
The first terms agree, because, because in Gekeler's notation, $t_p=t(pz)$, which translates to $u_p = u(p\omega_1)$ in our notation. Also remember that Gekeler's definition of the Hecke operator differs from ours by a scalar factor of $p^k$.
The second terms agree, but look different, because we normalize by a power of $p$ at a different stage in the computation. For example, in our computations, $L_b= e_{pA}(A)$, while Gekeler uses the finite lattice $\ker \rho_{pA}$, which is $e_{A}(p^{-1}A) = \frac{1}{p}L_b$.
\end{Ex}

In some cases, the Goss polynomials $G_{n,L_b}$ do not depend on the lattice $L_b$. For example, if $n\le q$, then $G_{n,L_b}(X)=X^n$. In such cases, the computation can be simplified into that of a lower rank Hecke operator. The term that simplifies is
\[ \sum_{\beta_{m,b}} f_n|_{k-n}\beta_{m,b}(\omega') G_{n}(u) = T_{k-n,\delta'}f(\omega')G_n(u),\]
where $T_{k-n,\delta'}$ is the weight $k-n$ Hecke operator associated with the double coset $\Gamma'\delta'\Gamma'$, where $\delta'$ is the $(r-1)\times(r-1)$ diagonal matrix with determinant $p$.

\begin{Cor}
\label{cor:low_exp}
The constant and linear terms of the $u$-expansion of $T_{\delta}f$ can be computed as
\[ T_{k,\delta}f = T_{k,\delta'}f_0 + T_{k-1,\delta'}f_1u + \cdots.\]
\end{Cor}
\begin{Proof}
By the remark preceding this Corollary, these two terms can be found from the computation in Theorem \ref{Thm:HeckeExp} in the second summand when $n=0$ and $n=1$. However, by Proposition \ref{Prop:GossPolynomials} (7), if $n>2$, there are no $X^0$ or $X^1$ terms in the Goss polynomial $G_{n,b}(X)$, so these are the only contributions to those terms.
\end{Proof}

%
%\begin{Rem}
%One might hope that Theorem \ref{Thm:HeckeExp} gives a way to inductively compute the action of these Hecke operators on Drinfeld modular forms for $\GL_r(\mathbb{F}_q[t])$. In general, the coefficients $f_n$ are only weak modular forms, so a similar result will be needed for weak modular forms that have a $u$-expansion that is not a power series. As far as I can tell, the only step that doesn't carry over is the step where we show that the sum $\sum_{b_1}u_{\Lambda'_b}(\omega_1+b_1)^n = G_n(u)$, because Goss polynomials only work for sums of the kind $\sum (z+\lambda)^{-n}$, where $n>0$. However, all may not be lost, since 
%\begin{align*} 
  %\sum_{b_1}e_{\Lambda'_b}(\omega_1+b_1\omega_m)^n &= \sum_{b_1}\left(e_{\Lambda'_b}(\omega_1)+e_{\Lambda'_b}(b_1\omega_m)\right)^n\\
	%&= \sum_{i=0}^n \binom{n}{i}e_{\Lambda'_b}(\omega_1)^{n-i}\sum_{b_1}e_{\Lambda'_b}(b_1\omega_m)^i,
	%\end{align*}
%by the Binomial Theorem and it may be possible to say something about this inner sum. 
%\end{Rem}

\begin{Cor}
\begin{enumerate}
	\item $T_\delta$ maps modular forms to modular forms.
	\item $T_\delta$ maps cusp forms to cusp forms.
	\item $T_\delta$ maps double cusp forms to double cusp forms.
\end{enumerate}
\end{Cor}

\section{Eigenforms and their eigenvalues}
\label{sec:eigen}

%\begin{Prop}\label{prop:EisEigen}
%Let $p\in A$ be irreducible. For every positive integer $k$ divisible by $q-1$, the Eisenstein series $E^k$ is an eigenform for $T_p$ with eigenvalue $p^k$. 
%\end{Prop}
%\begin{Proof}
%We use the description from Proposition \ref{Prop:HeckeActionComputation}. We claim that calculating $T_pE^k(\omega)$ is the same as evaluating
%\[\sum_{\lambda\in \frac{1}{p}\omega\Lambda}\lambda^{-k}\]
%which is clearly what we want. To show this we count the number of times each such $\lambda$ occurs in the sums
%\[E^k(\beta\omega)=\sumprm_{(a_1,\ldots,a_r)}(a_1\omega_1+\cdots+a_m\omega_{m,\beta}+\cdots+a_r)^{-k}.\]
%Note that the term 
%\[
%p^kE^k((p\omega_1,\ldots,p\omega_{r-1},1))=\sum_{\lambda\in \omega_1A+\cdots+\omega_{r-1}A+\frac{1}{p}A} \lambda^{-k},
%\]
%so each lattice $L\subset \Cinf$ with $L\supset \Lambda$ and $L/\Lambda\cong (A/pA)$ is counted exactly once.
%
%If $\lambda\in\Lambda$, then $\lambda$ is counted once in each of the $1+q^{\deg p}+\cdots+q^{(r-1)\deg p}$ times. Since we are in characteristic $p$, this contributes exactly $\lambda^{-k}$. 
%
%If $\lambda=\frac{a_1\omega_1+\cdots+a_r}{p}\notin \Lambda$, then set $m=\min\{i:a_i\ne 0\}$. Then there is a unique $\BF_q$ multiple of $\lambda$ of the form $\frac{\omega_m-b_{m+1}\omega_{m+1}-\cdots-b_r}{p}$ and $\beta$, and hence the superlattice $L$, can be uniquely determined from it. Therefore, each such $\lambda$ is counted exactly once. 
%\end{Proof}

\begin{Thm}\label{thm:DiscEigen}
Let $p\in \BF_q[t]$. The discriminant function $\Delta:\Omega^r\to\Cinf$ is an eigenform for $T_{p}$ with eigenvalue $p^{q^{r-2}(q-1)}$.
\end{Thm}
\begin{Proof}
By \cite[Cor 17.10]{BBP}, the space of cusp forms of weight $q^r-1$ is one-dimensional, so $\Delta$ must be an eigenform. It remains to compute its eigenvalue. 

Suppose that $\Delta$ has the $u$-expansion
\[ \Delta(\omega) = \sum_{n\ge q-1}f_n(\omega')u^{n}.\]
By the product formula for $\Delta$ (\cite[Thm 8]{BassonProd}) 
\[ \Delta(\omega) = -\Delta'(\omega')^qu^{q-1}\prod_{a\in A_+}(1+h_a(u))^{(q^r-1)(q-1)}\]
we know that $f_{q-1} = -\Delta'(\omega')^q$.  By Theorem \ref{Thm:HeckeExp} we get
\[ T_{\delta,q^r-1}\Delta(\omega) = \sum_{n} T_{\delta',q^r-1-n}f_n(\omega')G_n(u) + \sum_n f_n(\omega')u(p\omega_1)^n.\]
First we note that the second summand cannot contribute to the coefficient of $u^{q-1}$, since $u(p\omega_1)$ has degree $q^{\deg p - 1}>1$ in $u$ and $n\ge q-1$. For the first summand, we show that only the term $n=q-1$ contributes.

If $n=q-1$, then by Proposition \ref{Prop:GossPolynomials}(\ref{GP:SmallPowers}) we have $G_{q-1}(X)=X^{q-1}$  and 
%hence contributes the term $-p^{q-1}\Delta'(\omega')u^{q-1}$. To complete the proof we show that for
we claim that if $n\ge q$, then the polynomial $G_{n}(X)$ has no $X^{q-1}$ term. By Proposition \ref{Prop:CoefficientsCongruent}, the only values of $n$ for which $f_n\ne 0$ in the expansion of $\Delta$ satisfy $n\equiv 0$ or $n\equiv -1$ modulo $q$. If $q$ divides $n$, then every exponent in $G_{n}(X)$ is divisible by $q$, and the smallest such term is $X^q$. 

So suppose that $n=mq-1$ where $m\ge 2$. Then by Proposition \ref{Prop:GossPolynomials}(\ref{GP:Derivative}) we have
\[
X^2G'_{n}(X)=(q-1)G_{mq}(X)=-(G_m(X))^q.
\]
By Proposition \ref{Prop:GossPolynomials}(\ref{GP:NoXTerm}), if $m\ge 2$, then $G_m(X)$ has no $X$ term, and so the smallest possible term in $P_m(X)^q$ is $X^{2q}$ and hence $G_n(X)$ can have no $X^{q-1}$ term.

Thus, the leading term in the expansion of $T_{\delta, q^r-1}\Delta$ is $T_{\delta', q^r-q}\left(-\Delta'(\omega')^q\right)u^{q-1}$. 

\medskip

For the remainder of the proof, we proceed by induction on the rank. Gekeler \cite[Cor. 7.5]{GekelerCoeff} showed that in the rank 2 case, $\Delta$ is an eigenform with eigenvalue $p^{q-1}$. If the statement is assumed to be true in the rank $r-1$ case, then $\Delta'$ is an eigenform with eigenvalue $p^{q^{r-3}(q-1)}$, so
\[ T_{\delta', q^{r-1}-1}\Delta' = p^{q^{r-3}(q-1)}\Delta'.\]
Then 
\begin{align*} 
  T_{\delta',q^r-q}(\Delta')^q(\omega') &= \sum_{\beta} (\Delta')^q|_{q^r-q}\beta (\omega')\\
&= \sum_{\beta} (\Delta')^q(\beta(\omega'))j(\beta,\omega')^{-(q^r-q)}\\
& = \left(\sum_{\beta} \Delta'(\beta(\omega'))j(\beta,\omega')^{-(q^{r-1}-1)}\right)^q\\
&= \left( T_{\delta',q^{r-1}-1}\Delta'(\omega')\right)^q\\
&= (p^{q^{r-3}(q-1)}\Delta')^q = p^{q^{r-2}(q-1)}
\end{align*}
which finishes the proof by induction.
\end{Proof}

\begin{Thm}\label{thm:CoeffEigen}
The coefficient forms $g_1,\ldots,g_r$ are all eigenforms for $T_\delta$ and 
\begin{enumerate}
	\item the eigenvalue of $g_i$ is $p^{q^{i-2}(q-1)}$ if $i\ge 2$, 
	\item the eigenvalue of $g_1$ is $p^{q-1}$.
\end{enumerate}
\end{Thm}
\begin{Proof}
The coefficient form $g_1=(t^q-t)^{-1}E^{q-1}$ is a scalar multiple of an Eisenstein series, and therefore an eigenform. For $i=2,\ldots r-1$ we proceed by induction on the rank, and for $i=r$, we apply Theorem \ref{thm:DiscEigen}. Let $i$ be a fixed number from 2 to $r-1$. By the induction hypothesis $T_\delta g_i'=p^{q^{i-2}(q-1)}g_i'$. By \ref{cor:low_exp}, we know that
\[T_\delta g_i=(T_{\delta'}g_i')+O(u) = p^{q^{i-2}(q-1)}g_i' + O(u),\]
so that
\[T_\delta g_i-p^{q^{i-2}(q-1)}g_i\]
is a cusp form. But the lowest weight of any non-zero cusp form is $q^r-1>q^{i-1}$. Therefore $T_\delta g_i=p^{q-1}g_i$, and $g_i$ is an eigenform with eigenvalue $p^{q^{i-2}(q-1)}$. 
\end{Proof}

This idea of proving that a low weight modular form is an cusp form can also be used to prove that certain modular forms are not eigenforms.

\begin{Prop}\label{prop:NonExample}
Suppose that the characteristic of $\mathbb{F}_q$ is not equal to 2.
Let $r\ge 2$ and let $\delta$ be the $r\times r$ matrix with $t$ in the upper-left and the identity otherwise. Then the Drinfeld modular form $g_2^2$ is not an eigenform for $T_\delta$.
\end{Prop}
\begin{Proof}
The space of weight $2(q^2-1)$ forms is 3-dimensional for any rank, and spanned by the forms $g_2^2$, $g_1^{q+1}g_2$, and $g_1^{2q+2}$. Let us first show the statement in the rank 2 case. In that case, the forms $g_2^2$ and $g_1^{q+1}g_2$ are cusp forms, while $g_1^{2q-2}$ is not. Therefore $T_\delta(g_2^2)$ is a linear combination of $g_2^2$ and $g_1^{q+1}g_2$. Using the expansions in \cite[(6.4) and (10.3)]{GekelerCoeff} 
\[ g_1 = \bar\pi^{q-1}\left(1 -(t^q-t)u^{q-1}+ O(u^{q^2-q+1})\right),\qquad g_2 =  -\bar\pi^{q^2-1}\left(u^{q-1}-u^{q^2-q}+O(u^{q^2-q+1})\right),\]
where $\bar\pi$ is the Carlitz period, we see that their expansions start with 
\[g_2^2=\bar\pi^{2q^2-2}\left(u^{2q-2}-2u^{q^2-1}+\cdots\right)\quad\text{and}\quad g_1^{q+1}g_2=\bar\pi^{2q^2-2}\left(-1u^{q-1}+(t^q-t)u^{2q-2}+\cdots\right),\]
respectively. Thus, it is enough to show that the coefficient of $u^{q-1}$ in $T_t(g_2^2)$ is not zero.
%For anyone reading the comments, here are longer expansions
%$$g_1=1-(t^q-t)u^{q-1}-(t^q-t)u^{q^3-2q^2+2q-1}+\cdots $$
%$$g_2=-u^{q-1}+u^{q^2-q}-[1]u^{q^2-1}+\cdots $$
%$$g_1^{2q+2}=1-2(t^q-t)u^{q-1}+(t^q-t)^2u^{2q-2}-2(t^q-t)u^{q^2-q}+(t^q-t)u^{q^2-1}\cdots $$
%$$g_1^{q+1}g_2=-u^{q-1}+(t^q-t)u^{2q-2}+u^{q^2-q}+(t^q-t)u^{q^2-1}+\cdots $$
%$$g_2^2=u^{2q-2}-2u^{q^2-1}+2(t^q-t)u^{q^2+q}+\cdots $$
%and 
%$$T_pg_2^2= $$
For this we use Theorem \ref{Thm:HeckeExp} giving
\[ T_{\delta}\left(\sum a_iu^i\right) = \sum_{n} a_n u(t\omega')^n + \sum_n \sum_{b\in\mathbb{F}_q} a_nG_{n,L_b}(u),\]
where $L_b$ is the vector space $e_{tA}(A)$, one-dimensional over $\mathbb{F}_q$.

%which states that if $g_2^2$ has the expansion $\sum a_iu^i$, then there is the formula
%\[T_\delta\left(\sum a_iu^i\right)=p^{2q-2}\sum a_i u_t^i+\sum a_i G_{i,t}(tu)\]
%from. 

In this formula, $u(t\omega')$ is the parameter we defined in \eqref{eq:u_a}, and thus the first non-zero term in $\sum a_iu_t^i$ is $u^{2q-2}$ or higher. By \cite[(3.9)]{GekelerCoeff}, and since $t$ has degree 1, $G_{i}(X)$ is divisible by $X^q$ whenever $q\le \lfloor iq^{-1}\rfloor+1$, or equivalently $i\ge q^2-q$. By considering the $u$-expansion of $g_2^2$ we see that the only terms that influence the coefficient of $u^{q-1}$ are $a_{2q-2}G_{2q-2,L_b}(u)$. 

Since $L_b$ is one-dimensional, we can use \cite[(3.7)]{GekelerCoeff} to obtain the formula $G_{2q-2,t}(X)=X^{2q-2}-2\alpha X^{q-1}$, where $\alpha=-\lambda^{1-q}$ and $\lambda$ is the generator of the $\BF_q$-vector space $L_b$. 

We conclude that the coefficient of $u^{q-1}$ in $T_t(g_2^2)$ is $-2\alpha\ne 0$, and therefore $g_2^2$ is not an eigenform for $T_t$ in rank 2.

\medskip

Now we move to the higher rank case. Suppose that in rank 2 we have $T_\delta(g_2^2)=c_1g_1^{q+1}g_2+c_2g_2^{2}$, where $c_1,c_2\in \Cinf$ with $c_1\ne 0$ by the computation in the previous paragraph. Similarly to the previous Proposition we now prove by induction on $r$ that $g_2^2$ has this relation in any rank. Suppose that this is the case in rank $r-1$. In rank $r$, the constant coefficient in the expansion of $g_2^2$ is $g_2'^2$, while by Corollary \ref{cor:low_exp} that of $T_\delta(g_2^2)$ is $T_\delta(g_2'^2)=c_1g_1'^{q+1}g_2'+c_2g_2'^{2}$, which is the same as the constant coefficient of $c_1g_1^{q+1}g_2+c_2g_2^{2}$. Thus, $T_\delta(g_2^2)-(c_1g_1^{q+1}g_2+c_2g_2^{2})$ is a cusp form of weight $2q^2-2$, and therefore equals 0 if $r\ge 3$; and the same relation holds in rank $r$.
\end{Proof}

\section{Structure of the Hecke algebra}
\label{sct:HeckeAlgebra}
In \cite{ShimAut}, Shimura gives a very general setup for defining Hecke operators and studying Hecke algebras. We shall specialise his work from \cite[Section 3.1]{ShimAut} to our context and provide an analogue to his exposition of the case $\GL_r(\mathbb{Q})$ in \cite[Section 3.2]{ShimAut}. In this section we do not restrict $\delta$ to have only on non-trivial entry, but give only the first steps toward developing the theory for this more general case. In particular, we do not compute the action on Drinfeld modular forms for these more general $\delta$.

Let $\Gamma=\GL_r(A)$ and let $\Delta$ be the semi-group of $r\times r$ matrices with entries in $A$ and non-zero determinant. 
%Consider the double coset $\Gamma\delta\Gamma$, where $\delta\in \Delta$. (For now we do not assume that $\delta$ has the form \eqref{eq:DeltaForm}.) 
We define the Hecke algebra 
\[ R_\Gamma := \left\{\sum_{i} c_i\cdot \Gamma \delta\Gamma\mid c_i\in \mathbb{Z},\ \delta\in\Delta\right\}\]
as the free abelian group on the double cosets $\Gamma\alpha_i\Gamma$ for all $\alpha_i\in\Delta$. There is a multiplication operation on $R_\Gamma$. To define it, we write $\Gamma\alpha\Gamma=\cup \Gamma\alpha_i$ and $\Gamma\beta\Gamma = \cup \Gamma\beta_j$ and define
\UseTheoremCounterForNextEquation
\begin{equation}
\label{eq:DoubleCosetMult}
 (\Gamma \alpha\Gamma)\times (\Gamma\beta\Gamma) = \sum_\xi m(\alpha, \beta, \xi)\cdot(\Gamma\xi\Gamma),
\end{equation}
where $m(\alpha,\beta,\gamma)$ is the number of pairs $(i,j)$ for which $\Gamma\alpha_i\beta_j = \Gamma\xi$ and then extending $\mathbb{Z}$-linearly. Since the matrix transpose is an anti-automorphism of $\Gamma$, \cite[Prop 3.8]{ShimAut} tells us that the Hecke algebra is a commutative ring.

To aid our study, we shall make use of lattices inside the vector space $F^r$. In this section a lattice is a projective $A$-module of rank $r$ and a subset of $F^r$. The Theorem on Elementary Divisors (e.g. \cite[III. Thm 7.8]{LangAlgebra}) says that, because $A$ is a principal ideal domain, if we have two free modules $M\subset L$ of the same rank, then there exists a basis $\{e_1,\ldots, e_r\}$ of $L$ such that $M = a_1e_1A+\cdots a_re_rA$ for some $a_1,\ldots, a_r\in A$ such that $a_1$ is divisible by $a_2$, which is divisible by $a_3$, etc. These $a_i$ are uniquely determined up to units, so we can define an $A$-index of $M$ in $L$:

\begin{Def}
With the notation above, we define
\[ [L:M]_A := (a_1,a_2,\ldots, a_r).\]
\end{Def}
In that case $L/M \cong (A/a_1A)\times\cdots\times (A/a_rA)$. A consequence of the Theorem on Elementary Divisors is that one possible set of representatives for the double coset $\Gamma\backslash \Delta/\Gamma$ is the set of diagonal matrices $\diag(a_1,\ldots, a_r)$ with the divisibility property $a_r\mid a_{r-1}\mid\cdots \mid a_2\mid a_1$ (see \cite[Lemma 4.2.3]{BassonThesis} for more details). We shall often write $\Gamma\alpha\Gamma$ as $T(a_1,\ldots, a_r)$ to emphasise that we shall eventually think of it as an operator.

It is not hard to see that sublattices $M,N\subset A^r$ satisfy $[A^r:M]_A = [A^r:N]_A$ if and only if $M\alpha=N$ for some $\alpha\in \GL_r(A)$ (\cite[Lemma 4.2.5]{BassonThesis}, cf.~\cite[Lemma 3.12]{ShimAut}). A consequence is that if $\delta=\diag[a_1,\ldots, a_r]$ and if we write $\Gamma\delta\Gamma = \cup \Gamma\beta_j$, then the association $\Gamma\beta_j\to A^r\beta_j$ defines a bijection between the right cosets in $\Gamma\delta\Gamma$ and the lattices $M\subset A^r$ satisfying $[A^r:M]_A=(a_1,\ldots, a_r)$ (\cite[Lemma 4.2.6]{BassonThesis}, cf.~\cite[Lemma 3.13]{ShimAut}). 
%The number of right cosets in $\Gamma\delta\Gamma$ is called the degree of $\Gamma\delta\Gamma$, and thus equals the number of lattices $M\subset A^r$ of index $[A^r:M_A] = (a_1,\ldots, a_r)$. The degree map can be extended $\mathbb{Z}$-linearly to define a map $\deg:R_\Gamma \to \mathbb{Z}$ satisfying $\deg(x\cdot y) = \deg(x)\cdot \deg(y)$ (\cite[Prop 3.3]{ShimAut}).

When calculating the product of $\Gamma\alpha\Gamma=\cup \Gamma\alpha_i$ with $\Gamma\beta\Gamma=\cup \Gamma\beta_j$ in $R_\Gamma$, one can use the formula \eqref{eq:DoubleCosetMult}, but one needs to count $m(\alpha,\beta,\xi)$, the number of pairs $(i,j)$ for which $\Gamma\alpha_i\beta_j=\Gamma\xi$. 

\begin{Prop}
\label{prop:index}
 (\cite[Prop 4.2.8]{BassonThesis}, cf.~\cite[Prop 3.15]{ShimAut})
With $\alpha, \beta \in \Delta$, the number $m(\alpha, \beta,\xi)$ equals the number of lattices $M$ such that $[A^r:M]_A = [A^r:A^r\beta]_A$ and $[M:A^r\xi]_A = [A^r:A^r\alpha]_A$.
\end{Prop}
\begin{Proof}
With the lattice terminology, we may count instead the number of $(i,j)$ for which $A^r\alpha_i\beta_j=A^r\xi$, which is equivalent to $A^r\alpha_i = A^r\xi\beta_j^{-1}$.

First, we show that each pair $(i,j)$ that occurs in the count gives rise to a lattice $M$. Suppose $A^r\alpha_i\beta_j=A^r\xi$ and set $M=A^r\beta_j$. Then 
\[[A^r:M]_A = [A^r:A^r\beta_j]_A = [A^r:A^r\beta]_A\]
and
\[ [M:A^r\xi]_A = [A^r\beta_j:A^r\alpha_a\beta_j]_A = [A^r:A^r\alpha_i]_A = [A^r:A^r\alpha]_A.\]
Now we show that converse, that if $M$ is a lattice satisfying the formulas in the statement, then there exist a corresponding $i,j$ that contribute to the count. Since $[A^r:M]_A=[A^r:A^r\beta]$, the lattice $M$ has the same $A$-index in $A^r$ as $A^r\beta$ and must thus be one of the lattices $A^r\beta_j$ that corresponds to a coset $\Gamma\beta_j$ in the double coset $\Gamma\beta\Gamma$. This $\beta_j$ is uniquely determined and we have $[A^r:A^r\xi\beta_j^{-1}]_A = [A^r\beta_j:A^r\xi]_A=[A^r:A^r\alpha]_A$. Similarly $A^r\xi\beta_j^{-1}$ must now be a lattice corrresponding to some coset $\Gamma\alpha_i$ in the double coset $\Gamma\alpha\Gamma$, so $A^r\xi\beta_j^{-1}=A^r\alpha_i$. Thus $A^r\xi = A^r\alpha_i\beta_j$.
\end{Proof}

\begin{Prop}
\label{prop:Hecke_mult}
(\cite[Prop 4.2.9]{BassonThesis}, cf.~\cite[Prop 3.16]{ShimAut})
If $\det \alpha$ and $\det\beta$ are relatively prime in $A$, then $\Gamma\alpha\Gamma\beta\Gamma = \Gamma\alpha\beta\Gamma$.
%, or equivalently $T(a_1,\ldots, a_r)T(b_1,\ldots, b_r) = T(a_1b_1,\ldots, a_rb_r)$.
\end{Prop}
\begin{Proof}
Suppose that $[A^r:A^r\alpha]_A = (a_1,\ldots, a_r)$ and $[A^r:A^r\beta]_A = (b_1, \ldots, b_r)$. Then $\det\alpha = a_1a_2\cdots a_r$ and $\det\beta = b_1b_2\cdots b_r$, so the assumption implies that $a_1a_2\cdots a_r$ is relatively prime to $b_1b_2\cdots b_r$. The result follows if we can show that for any $\xi\in \Gamma\alpha\Gamma\beta\Gamma$, the multiplicity $m(\alpha, \beta,\xi)$ equals 1. By Proposition \ref{prop:index}, we need only show that for each $\xi$ there is a unique $M$ satisfying $[A^r:M]_A = (b_1,\ldots, b_r)$ and $[M:A^r\xi]_A = (a_1,\ldots, a_r)$. But such an $M$ can be uniquely constructed as $\frac{1}{a_1}A^r\xi\cap A^r$,  because this removes all $a_1$-torsion from the quotient $A^r/A^r\xi$, while not influencing the $b_1$-torsion at all, because $\gcd(a_1,b_1)=1$.
\end{Proof}

This means that we may focus our attention on double cosets of the form $\Gamma\diag[p^{e_1},\ldots, p^{e_r}]\Gamma$, where $p\in A$ is irreducible. Abbreviate this element as $T(p^{e_1},\ldots, p^{e_r})$.
%From the trivial relation $\Gamma (c\cdot I)\Gamma\beta\Gamma = \Gamma c\beta \Gamma$, we may even assume that $e_r=1$.
Define $R_{p}$ to be the subalgebra of $R_\Gamma$ generated by such double cosets. To study this algebra, we distinguish the following: Define $T_i$ to be the element in the $R_p$ corresponding to the double coset $\Gamma\diag[p,\ldots,p,1,\ldots,1]\Gamma$, where the diagonal matrix has $i$ $p$'s and $r-i$ ones. In fact, we often emphasise the rank by writing $T^{(r)}_i$ (and similarly $R_p^{(r)}$ for the algebra), because we shall use an inductive argument to establish that $T_1, T_2, \ldots, T_r$ generate $R_p$. Define
\[ \Psi:R_p^{(r+1)}\to R_p^{(r)},\begin{cases} T(p^{e_1},\ldots, p^{e_r},1)\mapsto  T(p^{e_1},\ldots, p^{e_r})&\\ T(p^{e_1},\ldots, p^{e_{r+1}}) \mapsto 0&\text{if }e_{r+1}>0.\end{cases}\]
and extend $\mathbb{Z}$-linearly to a $\mathbb{Z}$-module homomorphism.

\begin{Prop}
\label{prop:PsiHom}
(\cite[Lemma 4.2.11]{BassonThesis}, cf.~\cite[Lemma 3.19]{ShimAut})
The map $\Psi$ is a surjective ring homomorphism with kernel generated by $T^{(r+1)}_{r+1}$.
\end{Prop}
\begin{Proof}
Surjectivity is trivial and the shape of the kernel follows from the trivial relation $T(p,\ldots, p)T(a_1,\ldots, a_r) = T(pa_1,\ldots,pa_r)$.

It remains to show multiplicativity, i.e. to show that if $a=T(p^{a_1},\ldots, p^{a_r},1)$ and $b=T(p^{b_1},\ldots, p^{b_r},1)$, then $\Psi(a\cdot b) = \Psi(a)\cdot \Psi(b)$. To do this, we need to show that for every $c=T(p^{c_1},\ldots, p^{c_r},1)$, the multiplicities $m(a,b,c)$ and $m(\Psi(a),\Psi(b),\Psi(c))$ are equal. (Note that the product $a\cdot b$ also contains other terms, but they all have a last entry $p$ or larger, so map to 0 under $\Psi$.)

Choose $\xi\in\Delta$ such that $[A^{r+1}:A^{r+1}\xi]_A=c$. 
The multiplicity $m(a,b,c)$ is the number of lattices $M$ satisfying $[A^{r+1}:M]_A=b$ and $[M:A^{r+1}\xi]_A = a$. We assume that $(\{0\}^r\times A^r)\subset M$, since otherwise $c$ would have an $(r+1)$-st entry divisible by $p$, so maps to 0 under $\Psi$. Define $M'$ as the projection of $M$ onto its first $r$ coordinates and $\xi'$ as the upper-left $r\times r$ submatrix of $\xi$. Then $[A^r:M']_A = \Psi(b)$ and $[M':A^r\xi']_a=\Psi(a)$. Since one can construct $M$ from $M'$ as $M=M'\times A$, this defines a one-to-one correspondence between lattices $M$ that occur in the count of $m(a,b,c)$ and lattices $M'$ that occur in the count of $m(\Psi(a),\Psi(b),\Psi(c))$.
\end{Proof}

\begin{Prop}
(\cite[Thm 4.2.13]{BassonThesis}, cf.~\cite[Thm 3.20]{ShimAut})
The ring $R_p^{(r)}$ is the polynomial ring in $r$ algebraically independent elements $T^{(r)}_i$ for $i=1,2,\ldots, r$.
\end{Prop}
\begin{Proof}
We proceed by induction on $r$. The case $r=1$ follows, because $T(p^e) = T(p)^e$. Now assume that the result in true for $r$ and consider some element $x\in R^{(r+1)}_p$. By the induction hypothesis, $\Psi(x)$ is a polynomial in the $T^{(r)}_i$ and thus by Prop \ref{prop:PsiHom}, $x$ is that same polynomial (but with $T^{(r+1)}_i$ instead of $T^{(r)}_i$) plus a multiple of $T^{(r+1)}_{r+1}=T(p,p,\ldots, p)$. The factor that is multiplied by $T^{(r+1)}_{r+1}$ is again an element in $R_p^{(r+1)}$, so the same argument can be applied to it. The process will eventually end, because the determinant of $\alpha$ in any term $\Gamma\alpha\Gamma$ is reduced by a factor of $p^{r+1}$ when factoring out $T^{(r+1)}_{r+1}$.

If there were some polynomial relation
\[ \Phi\left(T_1^{(r+1)},\ldots, T_{r+1}^{(r+1)}\right) = \sum_{d=m}^n \left(T^{(r+1)}_{r+1}\right)^d \Phi_d\left(T_1^{(r+1)},\ldots, T_r^{(r+1)}\right) = 0,\]
with $\Phi_m,\Phi_n\ne 0$, then we can divide by $\left(T^{(r+1)}_{r+1}\right)^d$, because $T^{(r+1)}_{r+1}$ is not a zero divisor and apply $\Psi$ to get the following relation in rank $r$:
\[ \Phi_m\left(T_1^{(r)},\ldots, T_r^{(r)}\right)=0.\]
By the induction hypothesis, this implies that $\Phi_m$ must be the zero polynomial, contradicting the way we set up the polynomial $\Phi$.
\end{Proof}

There is now the interesting question of how one can write a general element of $R_p$, say $T(p^{e_1},\ldots, p^{e_r})$ as a polynomial in the elements $T_i$. We shall not attempt to answer this, but instead focus on another set of operators that are closely analogous to the classical operators $T_p$.

\begin{Def}
We define 
\[ \mathcal{T}_{p^n} = \sum_{e_1+\cdots+e_r=n} T(p^{e_1},\ldots, p^{e_r}),\]
in other words, $\mathcal{T}_{p^n}$ is the sum of the double cosets $\Gamma\alpha\Gamma$, where $\det\alpha = p^n$.
\end{Def}

\begin{Rem}
Note that if $n=1$, then $\mathcal{T}_p = T(p,1,\ldots, 1)$, which reduces to the same element $T_\delta$ studied in sections \ref{sec:Hecke} and \ref{sec:eigen}.
\end{Rem}

Recall that since $p\in A$ is irreducible, the quotient $A/pA$ is a field. It is well-known that the number of $k$-dimensional subspaces of the $n$-dimensional vector space $(A/pA)^n$. Let $d=\deg p$ equals the $q$-binomial coefficient (or rather the $q^d$-binomial coefficient
\[\qBinom{n}{k} = \frac{(q^{nd}-1)(q^{nd}-q^d)\cdots (q^{nd}-q^{(k-1)d})}{(q^{kd}-1)(q^{kd}-q^d)\cdots (q^{kd}-q^{(k-1)d}}\]
and it is also well-known that $\qBinom{n}{n-k} = \qBinom{n}{k}$, like for Binomial coefficients.

\begin{Lem}
\label{lem:inversion}
\begin{enumerate}
	\item The degree of $T_i$ equals $\qBinom{r}{i}$.
	\item For $k>0$, one has
	\[ \sum_{i=0}^k (-1)^iq^{\tfrac{1}{2}i(i-1)d}\qBinom{k}{i} = 0.\]
\end{enumerate}
\end{Lem}
\begin{Proof}
The degree of $T_i$ equals the number of lattices $M\subset A^r$ with $A$-index $[A^r:M]_A=(p,\ldots, p,1,\ldots, 1)$, with $i$ $p$'s. Such a lattice satisfies $pA^r\subset M\subset A^r$ and thus $M/pA^r\cong (A/pA)^{r-i}$. This correspondence between lattices $M$ and subspaces $M/pA^r$ of $A^r/pA^r\cong (A/pA)^r$ is one-to-one.

For the second statement, define the polynomials $\Phi(x) = (x-1)(x-q^d)\cdots(x-q^{(k-1)d})$ and 
\[ \Psi(x) = \sum_{i=0}^{k-1}\frac{\Phi(x)}{\Phi'(q^{id})(x-q^{id})}.\]
Then $\Psi$ has degree less than or equal to $k-1$ and satisfies $\Psi(q^{jd})=1$ for the $k$ values $j=0,1,\ldots, k-1$, and thus $\Psi(x) = 1$ identically. Then $\Psi(q^{kd})=1$ from which the statement can be deduced.
\end{Proof}

It is simpler to study the composites $T_i\mathcal{T}_{p^j}$ before moving on to the composites $\mathcal{T}_{p^i}\mathcal{T}_{p^j}$.

\begin{Lem}
\label{lem:composite}
Making the conventions $\qBinom{0}{0}=1$ and $\qBinom{k}{i}=0$ if $i>k$, one gets
\[ T_i\mathcal{T}_{p^j} = \sum_{k=0}^r \qBinom{k}{i}\sum_{d_1\ge\cdots\ge d_k\ge 1} T(p^{d_1},\ldots,p^{d_k},1,\ldots, 1),\]
where the sum is taken over those $d_i$ such that $d_1+\cdots+d_k=i+j$.
\end{Lem}
\begin{Proof}
By Prop \ref{prop:index}, we can compute the coefficient of $T(p^{d_1},\ldots p^{d_k},1\ldots, 1)$ by fixing some $\xi\in\Delta$ with $[A^r:A^r\xi]_A = (p^{d_1},\ldots, p^{d_k},1,\ldots, 1)$ and counting the number of lattices $M$ satisfying $[A^r:M]_A=(p,\ldots, p,1,\ldots,1)$. Note that Prop \ref{prop:index} also specifies a condition on $[A^r:A^r\xi]_A$, but this condition is irrelevant, because by the way $\mathcal{T}_{p^j}$ is defined, all possible $A$-indices $[A^r:A^r\xi]_A$ are allowed. If $k\ge i$, this means that $pA^r+A^r\xi\subseteq M\subseteq A^r$. Since $A^r/(A^r\xi)\cong (A/pA)^{d_1}\times\cdots\times (A/pA)^{d_k}$, we deduce that $A^r/(pA^r+A^r\xi)\cong (A/pA)^k$. We also know that $A^r/M\cong (A/pA)^i$ and thus $M/(pA+A^r\xi)\cong (A/pA)^{k-i}$. Thus, the number of lattices $M$ is equal to the number of subspaces isomorphic to $(A/pA)^{k-i}$ in $(A/pA)^k$. If $k<i$, then it is impossible to find a suitable $M$, so the count is 0.
\end{Proof}

\begin{Prop}
\label{prop:relation}
For any $k\ge r$ we have
\[ \mathcal{T}_{p^k} = \sum_{i=1}^r (-1)^{i+1}q^{\tfrac{1}{2}i(i-1)d}T_i\mathcal{T}_{p^{k-i}}.\]
\end{Prop}
\begin{Proof}
One verifies this by direct substitution of the results from Lemmas \ref{lem:inversion} and \ref{lem:composite}:

\begin{align*}
\sum_{i=1}^r (-1)^{i+1}q^{\tfrac{1}{2}i(i-1)d}T_i\mathcal{T}_{p^{k-i}} &= \sum_{i=1}^r (-1)^{i+1}q^{\tfrac{1}{2}i(i-1)d}\sum_{j=0}^r\qBinom{j}{i}\sum T(p^{d_1},\ldots, p^{d_j},1,\ldots, 1)\\
&=\sum_{j=0}^r\sum_{i=1}^j (-1)^{i+1}q^{\tfrac{1}{2}i(i-1)d}\qBinom{j}{i} \sum T(p^{d_1},\ldots, p^{d_j},1,\ldots, 1)
\end{align*}
The sums $\sum_{i=0}^j (-1)^{i+1}q^{\tfrac{1}{2}i(i-1)d}\qBinom{j}{i}$ equal 0 unless $i=j=0$, so we remain with (negative) the terms corresponding to $i=0$:
\[ \sum_{j=0}^r\sum T(p^{d_1},\ldots, p^{d_j},1,\ldots, 1).\]
As $j$ ranges from 1 to $r$, we get all possible elements 
$T(p^{d_1},\ldots, p^{d_j},1,\ldots, 1)$ with $d_1+\cdots+d_j=k$, which by definition equals $\mathcal{T}_{p^k}$.
\end{Proof}

\begin{Thm}
\label{thm:HeckeMult}
The Hecke operators $\mathcal{T}_{p^k}$ act in a completely multiplicative way on Drinfeld modular forms.
\end{Thm}
\begin{Proof}
Since Drinfeld modular forms are functions with a codomain of finite characteristic, we need to consider the relation from Prop \ref{prop:relation} modulo the characteristic to get
\[ \mathcal{T}_{p^k} = T_i\mathcal{T}_{p^{k-i}}\quad \pmod{\text{char }{\mathbb{F}_q}}.\]
In particular, when $i=1$, then $T_1 = \mathcal{T}_p$ and thus $\mathcal{T}_{p^k} = \mathcal{T}_p\mathcal{T}_{p^{k-1}}$ and inductively it follows that $\mathcal{T}_{p^k} = (\mathcal{T}_p)^k$. 

The multiplicativity $\mathcal{T}_M\mathcal{T}_N = \mathcal{T}_{MN}$ in the case $\gcd(M,N)=1$ follows by applying Proposition \ref{prop:Hecke_mult} to each pair of terms from the definitions of $\mathcal{T}_M$ and $\mathcal{T}_N$.
\end{Proof}

%%%%%%%%%%%%%%%%%%%%%%%%%%%%%%%%%%%%%%%%%%%%%%%%%%%%%%%%%%%%%%%%%%%%%%%%%%

%%%%%%%%%%%%%%%%%%%%%%%%%%%%%%%%%%%%%%%%%%%%%%%%%%%%%%%%%%%%%%%%%%%%%%%%%%%%%%%%%%%%%%%%%%%%%%%%%%%%%%%%%%%%%%%%%%%%%%%%%%%%%%%%%%%

\begin{center}
\rule{8cm}{0.01cm}

\begin{minipage}[t]{6.5cm}{\small
Department of Mathematical Sciences \\
University of Stellenbosch \\
Stellenbosch, 7600 \\
South Africa \\
djbasson@sun.ac.za}
\end{minipage}\hfill
\end{center}

\end{document}